\def\titlep{Jump transformations and an embedding of ${\cal O}_{\infty}$ into ${\cal O}_{2}$}
\documentclass[11pt]{article}
\usepackage{graphicx,ifthen}
\usepackage{amssymb}

\font\germ=eufm10 at12pt

\def\goth#1{\hbox{\germ#1}}

\newcommand{\qed}{\hbox{\rule[-2pt]{3pt}{6pt}}}
\newcommand{\qedh}{\hfill\qed \\}

\newtheorem{Thm}{Theorem}[section]

\newtheorem{rem}[Thm]{Remark}

\newtheorem{defi}[Thm]{Definition}
\newtheorem{lem}[Thm]{Lemma}

\newtheorem{cor}[Thm]{Corollary}

\newtheorem{prob}[Thm]{Problem}

\def\cal#1{\mathcal #1}
\def\con{{\cal O}_{N}}
\def\coni{{\cal O}_{\infty}}
\def\edot{=1,\ldots,N}
\def\pr{{\it Proof.}\quad}

\def\co#1{{\cal O}_{#1}}
\def\lxm{L_{2}(X,\mu)}

\def\disp#1{{\displaystyle #1}}

\def\rnd{Radon-Nikod\'{y}m
derivative}

\def\bfs{branching function system}

\addtocounter{footnote}{1}
\def\sftt#1{
\setcounter{equation}{0}
\addtocounter{footnote}{1}
\section{#1}
}

\def\ssft#1{\subsection{#1}}
\def\sssft#1{\subsubsection{#1}}

\def\cls{\quad
\clearpage
}

\begin{document}
\def\inputcls#1{\cls\input #1.txt}
\def\plan#1#2{\par\noindent\makebox[.5in][c]{#1}
\makebox[.1in][l]{$|$}
\makebox[3in][l]{#2}\\}
\def\nset#1{\{1,\ldots,N\}^{#1}}
\def\brl{branching law}
\def\bfsnl{{\rm BFS}_{N}(\Lambda)}
\def\scm#1{S({\bf C}^{N})^{\otimes #1}}
\def\mqb{\{(M_{i},q_{i},B_{i})\}_{i=1}^{N}}
\newcommand{\mline}{\noindent
\thicklines
\setlength{\unitlength}{.18mm}
\begin{picture}(1000,5)
\put(0,0){\line(1,0){750}}
\end{picture}
}
\def\boxtimes{\noindent
\setlength{\unitlength}{.020918mm}
\begin{picture}(120,150)(60,0)
\thinlines
\put(0,0){
\line(1,0){100}\line(0,1){100}
}
\put(100,100){
\line(-1,0){100}\line(0,-1){100}
}
\put(40,0){$\times$}
\end{picture}
 }
\def\authorp{Katsunori  Kawamura}
\def\authorq{Dan Lascu}
\def\authorr{Ion Coltescu}
\def\emailp{Electronic mail: kawamura@kurims.kyoto-u.ac.jp.}
\def\emailq{Electronic mail: lascudan@gmail.com.}
\def\emailr{Electronic mail: icoltescu@yahoo.com.}
\def\addressp{{\small {\it College of Science and Engineering Ritsumeikan University,}}\\
{\small {\it 1-1-1 Noji Higashi, Kusatsu, Shiga 525-8577, Japan,}}
}
\def\addressq{{\small {\it 
Mircea cel Batran Naval Academy, 1 Fulgerului, 900218 Constanta, Romania}}}
\def\addressr{****}
\def\ba{\mbox{\boldmath$a$}}
\def\bb{\mbox{\boldmath$b$}}
\def\bc{\mbox{\boldmath$c$}}
\def\be{\mbox{\boldmath$e$}}
\def\bp{\mbox{\boldmath$p$}}
\def\bq{\mbox{\boldmath$q$}}
\def\bu{\mbox{\boldmath$u$}}
\def\bv{\mbox{\boldmath$v$}}
\def\bw{\mbox{\boldmath$w$}}
\def\bx{\mbox{\boldmath$x$}}
\def\by{\mbox{\boldmath$y$}}
\def\bz{\mbox{\boldmath$z$}}
\def\aei{almost everywhere in}

%
%
\pagestyle{plain}
\setcounter{page}{1}
\setcounter{section}{0}

\title{\titlep}
\author{
Katsunori Kawamura\thanks{\emailp}\\
\addressp\\
\\
Dan Lascu\thanks{\emailq}\, and 
Ion Coltescu\thanks{\emailr}\\
\addressq\\
}
\date{}
\maketitle
%
%
\begin{abstract}
A measurable map $T$ on a measure space induces 
a representation $\Pi_{T}$ of a Cuntz algebra ${\cal O}_{N}$
when $T$ satisfies a certain condition.
For such two maps $\tau$ and $\sigma$ and
representations $\Pi_{\tau}$ and $\Pi_{\sigma}$ associated with them,
we show that $\Pi_{\tau}$ is the restriction of $\Pi_{\sigma}$
when $\tau$ is a jump transformation of $\sigma$.
Especially, the Gauss map $\tau_1$ and the Farey map $\sigma_1$ 
induce representations 
$\Pi_{\tau_1}$ of ${\cal O}_{\infty}$ and that $\Pi_{\sigma_1}$ 
of ${\cal O}_{2}$, respectively,
and $\Pi_{\tau_1}=\Pi_{\sigma_1}|_{{\cal O}_{\infty}}$ with respect to
a certain embedding of ${\cal O}_{\infty}$ into ${\cal O}_{2}$.
\end{abstract}

\sftt{Introduction}
\label{section:first}
The purpose of this paper is to show a new relation 
between dynamical systems and operator algebras.
The former means jump transformations in metric number theory
and the later does an embedding of Cuntz algebras.
In this section, we show our motivation and main theorem.
%
%
\ssft{Motivation}
\label{subsection:firstone}
We explain our motivation in this subsection.
Mathematical details will be shown in $\S$ \ref{subsection:firsttwo},
$\S$ \ref{subsection:firstthree} and $\S$ \ref{section:third}.

As is well known, metric number theory has originated as Gauss' problem
about the asymptotic behavior of iterations of the regular continued fraction transformation.
As a modern style,
the theory is formulated as a measure theoretical dynamical system \cite{IK,Schweiger}.
For a measure space $(X,\mu)$, let $T$ be a measurable map from $X$ to $X$
which may not be invertible.
We call the triplet $(X,\mu,T)$ a {\it (measure theoretical) dynamical system}.
For example,
the {\it Gauss map (or the continued fraction transformation)} $\tau_{1}$ 
on the closed interval $[0,1]$ defined by
%
%
\begin{equation}
\label{eqn:gauss}
\tau_{1}(x)\equiv 
\left\{
\begin{array}{ll}
\disp{\frac{1}{x}-\left\lfloor \frac{1}{x}\right\rfloor}
\quad&(x\ne 0),\\
\\
0 \quad&(x=0),\\
\end{array}
\right.
\end{equation}
gives an important dynamical system $([0,1],\lambda,\tau_{1})$
with respect to the Lebesgue measure $\lambda$
where $\left\lfloor \cdot \right\rfloor$ denotes 
Gauss' symbol (or the floor function). 
The second and the third authors have studied
a generalization of the continued fraction transformation in metric number theory \cite{CL}.

On the other hand, the first author has studied representations of Cuntz algebras 
and Perron-Frobenius operators \cite{PFO01}.
A dynamical system $(X,\mu,T)$ induces 
a representation $(L_{2}(X,\mu),\Pi_{T})$ of a Cuntz algebra
when $(X,\mu,T)$ satisfies a certain condition.
Properties of $(L_{2}(X,\mu),\Pi_{T})$ are
characterized by $(X,\mu,T)$.
As an application,
such construction is used to construct type ${\rm III}$ factor representations
of Cuntz-Krieger algebras (\cite{TS11}, $\S$ 1.3).

These two different studies have the following common problem in dynamical system.
%
%
\begin{prob}
\label{prob:invariant}
For a given map $T$ on a set $X$,
find a $T$-invariant (or $T$-preserving) measure $\mu$ of $X$, that is,
%
%
\begin{equation}
\label{eqn:im}
\mu(T^{-1}(E))=\mu(E)
\end{equation}
for every $\mu$-measurable subset $E$ of $X$.
\end{prob}
For example,
a solution of Gauss' problem is given by an invariant measure of the Gauss map
\cite{IK,Schweiger},
and there exists a relation between an invariant measure
and the representation of Cuntz-Krieger algebra 
arising from dynamical system \cite{PFO01}.

The new idea in this paper is to show 
a relation between (the construction of) 
jump transformation and an embedding of $\coni$ into $\co{2}$
where $\coni$ and $\co{2}$ denote Cuntz algebras
which will be explained in $\S$ \ref{subsubsection:firsttwoone}.
Let $\tau$ and $\sigma$ be two maps, and
let $\Pi_{\tau}$ and $\Pi_{\sigma}$ denote representations associated with them,
respectively.
We show that $\Pi_{\tau}$ is the restriction of $\Pi_{\sigma}$
when $\tau$ is a jump transformation of $\sigma$ (Theorem \ref{Thm:maintwo}).
Especially, the Gauss map $\tau_1$ and the Farey map $\sigma_1$ 
induce representations 
$\Pi_{\tau_1}$ of ${\cal O}_{\infty}$ and that $\Pi_{\sigma_1}$ of ${\cal O}_{2}$, 
respectively,
and $\Pi_{\tau_1}=\Pi_{\sigma_1}|_{{\cal O}_{\infty}}$ with respect to
a certain embedding of ${\cal O}_{\infty}$ into ${\cal O}_{2}$
($\S$ \ref{subsection:thirdtwo}).

%
%
\ssft{Representations of Cuntz algebras arising from dynamical systems}
\label{subsection:firsttwo}
In this subsection, we explain representations of Cuntz algebras
arising from dynamical systems
according to \cite{PFO01}.
%
%
\sssft{Cuntz algebras}
\label{subsubsection:firsttwoone}
For $N=2,3,\ldots,+\infty$, 
let $\con$ denote the {\it Cuntz algebra} \cite{C}, that is, a C$^{*}$-algebra 
which is universally generated by $s_{1},\ldots,s_{N}$ satisfying
%
%
\begin{eqnarray}
\label{eqn:first}
s_{i}^{*}s_{j}=&\delta_{ij}I\quad &\mbox{ for }i,j\edot,\\
\label{eqn:seconda}
\sum_{i=1}^{N}s_{i}s_{i}^{*}=&I\quad&(\mbox{if } N<+\infty),\\
\label{eqn:secondb}
\sum_{i=1}^{k}s_{i}s_{i}^{*}\leq &I\quad &\mbox{ for }k= 1,2,\ldots\,(\mbox{if }N = +\infty)
\end{eqnarray}
where $I$ denotes the unit of $\con$.

Since $\con$ is simple, that is, there is no
nontrivial closed two-sided ideal,
any unital homomorphism from $\con$ to a C$^{*}$-algebra is injective.
If $t_{1},\ldots,t_{N}$ are elements of a unital C$^{*}$-algebra
${\goth A}$ such that
$t_{1},\ldots,t_{N}$ satisfy the relations of canonical generators of $\con$,
then the correspondence $s_{i}\mapsto t_{i}$ for $i\edot$
is uniquely extended to a $*$-embedding
of $\con$ into ${\goth A}$ from the uniqueness of $\con$.
Therefore we call such a correspondence 
among generators by an embedding of $\con$ into ${\goth A}$.

Assume that $s_{1},\ldots,s_{N}$ are realized as
operators on a Hilbert space ${\cal H}$.
According to (\ref{eqn:first}) and (\ref{eqn:seconda}),
${\cal H}$ is decomposed into orthogonal subspaces as
$s_{1}{\cal H}\oplus\cdots \oplus s_{N}{\cal H}$.
Since $s_{i}$ is an isometry,
$s_{i}{\cal H}$ has the same dimension as ${\cal H}$.
From this, we see that 
there is no finite dimensional representation of $\con$ which preserves 
the unit.
The following illustration is helpful in understanding $s_{1},\ldots,s_{N}$:

\def\firstbox{
\put(0,0){\line(1,0){500}}
\put(0,30){\line(1,0){500}}
\put(0,0){\line(0,1){30}}
\put(500,0){\line(0,1){30}}
\put(240,10){${\huge {\cal H}}$}
}
\def\secondbox{
\put(0,0){\line(1,0){500}}
\put(0,30){\line(1,0){500}}
\put(0,0){\line(0,1){30}}
\put(100,0){\line(0,1){30}}
\put(200,0){\line(0,1){30}}
\put(300,0){\line(0,1){30}}
\put(140,10){$\cdots$}
\put(340,10){$\cdots$}
\put(400,0){\line(0,1){30}}
\put(500,0){\line(0,1){30}}
\put(40,10){${\huge s_{1}{\cal H}}$}
\put(240,10){${\huge s_{i}{\cal H}}$}
\put(440,10){${\huge s_{N}{\cal H}}$}
}
\def\cross{
\thinlines
\put(225,60){$s_{i}$}
\put(245,60){$\downarrow$}
\qbezier[200](500,100)(400,65)(300,30)
\qbezier[200](0,100)(100,65)(200,30)
}
\setlength{\unitlength}{.22mm}
\begin{picture}(1000,150)(0,-10)
\thicklines
\put(0,100){\firstbox}
\put(0,0){\secondbox}
\put(0,0){\cross}
\end{picture}

%
%
\sssft{Representations of Cuntz algebras}
\label{subsubsection:firsttwotwo}
We introduce \bfs s.
Let $(X,\mu)$ be a measure space and 
$\Omega_{N}\equiv \{1,\ldots,N\}$ for $2\leq N<\infty$
and $\Omega_{\infty}\equiv {\bf N}=\{1,2,3,\ldots\}$.
Let $T$ be a measurable map on $(X,\mu)$. 
Define the new measure $\mu\circ T$ on $X$ by
$(\mu\circ T)(E)\equiv \mu(T(E))$ for $E\subset X$.
For $2\leq N\leq \infty$,
a family $f=\{f_{i}:i\in \Omega_{N}\}$ of maps on $X$
is a {\it \bfs} if 
\begin{enumerate}
\item
$f_{i}$ is a measurable map from $X$ to $X$ for each $i$,
\item
if $R_{i}\equiv f_{i}(X)$, then
$\mu\bigl(X\setminus \bigcup_{i\in \Omega_{N}} R_{i}\bigr)=0$ and
$\mu(R_{i}\cap R_{j})=0$
when $i\ne j$, and
\item
there exists the \rnd\ $\Phi_{f_{i}}$ of $\mu\circ f_{i}$
with respect to $\mu$ and $\Phi_{f_{i}}>0$ \aei\ $X$ for each $i$.
\end{enumerate}
A map $F$ on $X$ is called the
{\it coding map} of a \bfs\ $f=\{f_{i}:i\in \Omega_{N}\}$
if $F\circ f_{i}=id_{X}$ for each $i$.

Let $\lxm$ denote the Hilbert space
of all square integrable complex valued-functions on $X$.
For a \bfs\ $f=\{f_{i}:i\in \Omega_{N}\}$
with the coding map $F$,
define the family $\{S(f_{i}):i\in\Omega_{N}\}$
of operators on $\lxm$ by
%
%
\begin{equation}
\label{eqn:pfotwo}
\{S(f_{i})\phi\}(x)\equiv
\chi_{R_{i}}(x)\cdot 
\{\Phi_{F}(x)
\}^{1/2}\cdot 
\phi(F(x))\quad(\phi\in \lxm,\,x\in X)
\end{equation}
where $\chi_{R_{i}}$ denotes the characteristic function of $R_{i}$.
From this, 
adjoint operators $\{S(f_{i})^{*}:i\in\Omega_{N}\}$ are given as follows:
%
%
\begin{equation}
\label{eqn:conjugate}
\{S(f_{i})^{*}\phi\}(x)=\{\Phi_{f_{i}}(x)\}^{1/2}\cdot 
\phi(f_{i}(x))\quad(\phi\in \lxm,\,x\in X).
\end{equation}
Then $S(f_{i})$ is an isometry and
%
%
\begin{equation}
\label{eqn:homomorphism}
S(f_{i}) S(f_{j})=S(f_{i}\circ f_{j})\quad(i,j\in\Omega_{N})
\end{equation}
where $(f_{i}\circ f_{j})(x)\equiv f_{i}(f_{j}(x))$.
Let $\{s_{i}:i\in\Omega_{N}\}$ denote canonical generators of $\con$.
For a \bfs\ $f=\{f_{i}:i\in\Omega_{N}\}$,
%
%
\begin{equation}
\label{eqn:permurep}
\pi_{f}(s_{i})\equiv S(f_{i})\quad(i\in\Omega_{N}),
\end{equation}
defines a representation
$(\lxm,\pi_{f})$ of the Cuntz algebra $\con$
because
$\{S(f_{i}):i\in \Omega_{N}\}$ satisfy 
relations of canonical generators of $\con$.
Especially,
\bfs s $f=\{f_{1},f_{2}\}$ and 
$g=\{g_{i}:i\in{\bf N}\}$ 
define a representation of $\co{2}$ and 
that of $\coni$, respectively.
Examples will be shown in $\S$ \ref{section:third}.
%
%
\begin{rem}
\label{rem:one}
{\rm
For a dynamical system $(X,\mu,T)$,
a branching function system with the coding map $T$ is not unique
even if it exists.
Therefore 
the representation $\pi$ of $\con$ arising from $(X,\mu,T)$
depends on the choice of 
a branching function system with the coding map $T$.
}
\end{rem} 
%
%
\ssft{Main theorem}
\label{subsection:firstthree}
In this subsection, we show our main theorem.
For this purpose, we introduce jump transformation
according to \cite{Schweiger}.
%
%
\begin{defi}
\label{defi:jump}
Let  $(X,\mu,T)$ be a dynamical system.
Assume that $A$ is a measurable subset of $X$ such that 
%
%
\begin{equation}
\label{eqn:inverse}
\mu\biggl(X\setminus \bigcup_{n\geq 1}T^{-n}(A)\biggr)=0.
\end{equation}
\begin{enumerate}
\item
For the pair $(T,A)$,
the map $e$ from $X$ to ${\bf Z}_{\geq 0}\equiv\{0,1,2,3,\ldots\}$ 
is called the first entry time of $(T,A)$ if 
%
%
\begin{equation}
\label{eqn:time}
e(x)\equiv \min\{k\in {\bf Z}_{\geq 0}:T^{k}(x)\in A\}\quad(x\in X).
\end{equation}
%
\item
For the pair $(T,A)$,
define the map $J_{T,A}$ from $X$ to $X$ by
%
%
\begin{equation}
\label{eqn:jump}
J_{T,A}(x)\equiv T^{e(x)+1}(x)\quad(x\in X).
\end{equation}
We call $J_{T,A}$ the jump transformation of $(T,A)$. 
\end{enumerate}
\end{defi}
We briefly explain the meaning of jump transformation as follows.
The assumption (\ref{eqn:inverse}) means that
it surely takes finite time such that $x$ enters the subset $A$
with respect to the discrete time evolution
%
%
\begin{equation}
\label{eqn:tmapthree}
x\mapsto T(x)
\end{equation}
for almost all $x\in X$.
From this explanation, we can also understand the meaning of the first entry time.
For $x\in X$, assume $n\equiv e(x)<\infty$.
Then the transformation $J_{T,A}$ maps $x$ as follows:

\noindent
\def\jump{
\put(0,0){$x\mapsto T(x)\mapsto T^{2}(x)\mapsto\cdots\mapsto T^{n}(x)
\mapsto T^{n+1}(x)$}
\put(278,46){$A$}
\put(280,25){\rotatebox{90}{$\in$}}
\qbezier(20,-10)(200,-60)(360,-10)
\put(170,-60){$J_{T,A}$}
\put(365,-7){\vector(2,1){0}}
}
\thicklines
\setlength{\unitlength}{.18mm}
\begin{picture}(600,150)(-100,-75)
\put(0,0){\jump}
\end{picture}

Let $\{t_{1},t_{2}\}$ and  $\{s_{n}:n\in {\bf N}\}$ denote
canonical generators of Cuntz algebras $\co{2}$ and $\coni$, respectively.
Assume that $\coni$ is embedded into $\co{2}$ by
%
%
\begin{equation}
\label{eqn:embedding}
s_{n}= t_{2}^{n-1}t_{1}\quad(n\geq 1)
\end{equation}
where we write $t_{2}^{0}\equiv I$ for convenience.
This embedding induces the restriction of representations of 
$\co{2}$ on $\coni$:
\[{\rm Rep}\co{2}\ni\pi\mapsto \pi|_{\coni}\in{\rm Rep}\coni.\]
%
%
\begin{Thm}
\label{Thm:maintwo}
Let $(X,\mu,\sigma)$ be dynamical system
such that $\sigma$ is the coding map of some branching function system 
$f=\{f_{1},f_{2}\}$.
Assume that $A_{1}\equiv f_{1}(X)$ satisfies (\ref{eqn:inverse})
with respect to $(X,\mu,\sigma)$.
Let $\pi_{f}$ denote the representation
of $\co{2}$ associated with $f$ in (\ref{eqn:permurep}).
Let $J_{\sigma,A_{1}}$ be as in (\ref{eqn:jump}) for $(\sigma,A_{1})$.
Then the following holds:
\begin{enumerate}
\item
There exists a branching function system $g=\{g_{n}:n\in {\bf N}\}$ on $(X,\mu)$
with the coding map $J_{\sigma,A_{1}}$ such that
the representation $\pi_{g}$ of $\coni$ coincides
with the restriction $\pi_{f}|_{\coni}$ of $\pi_{f}$ on $\coni$
with respect to the embedding in (\ref{eqn:embedding}), that is,
%
%
\begin{equation}
\label{eqn:restriction}
\pi_{f}|_{\coni}=\pi_{g}.
\end{equation}
%
\item
If $\phi$ is the density of an invariant measure 
with respect to $\sigma$, 
that is,
%
%
\begin{equation}
\label{eqn:invariantone}
\nu(E)\equiv \int_{E}\phi(x)\,d\mu(x)\quad(E\subset X)
\end{equation}
defines a $\sigma$-invariant measure $\nu$ of $X$,
and $\phi$ belongs to $L_{1}(X,\mu)$,
then the function $\psi$ on $X$ defined by
%
%
\begin{equation}
\label{eqn:density}
\psi(x)=\left\{\{\pi_{f}(t_{1})^{*}\sqrt{\phi}\}(x)\right\}^{2}\quad(x\in X),
\end{equation}
is the density of an invariant measure of $J_{\sigma,A_{1}}$
where $\sqrt{\phi}(x)\equiv \sqrt{\phi(x)}$.
\end{enumerate}
\end{Thm}

We summarize Theorem \ref{Thm:maintwo}(i) as follows:
\[
\begin{array}{l|l}
\hline
\mbox{{\sf map}} &\mbox{{\sf Cuntz algebra}}\\
\hline
\sigma & \mbox{representation of }\co{2}=C^{*}\langle \{t_{1},t_{2}\}\rangle\\
\tau & \mbox{representation of }\coni=C^{*}\langle \{s_{n}:n\in {\bf N}\}\rangle\\
\tau=J_{\sigma,A_{1}}& s_{n}=t_{2}^{n-1}t_{1}\\
\hline
\end{array}
\]

%
%
\begin{rem}
{\rm
\begin{enumerate}
\item
We stress that jump transformation (\ref{eqn:jump})
and the embedding (\ref{eqn:embedding})
are independently introduced in different studies.
Theorem \ref{Thm:maintwo} is a casual discovery. 
We see that 
the index ``$e(x)+1$" of $T^{e(x)+1}(x)$
and the index ``$n-1$" of the element $t_{2}$ of the equation
$s_{n}=t_{2}^{n-1}t_{1}$ fit like a glove
in the proof of Theorem \ref{Thm:maintwo}.
This is an interesting accidental coincidence.
\item
From Theorem \ref{Thm:maintwo},
we see that a jump transformation gives
a geometric realization of the embedding 
of $\coni$ into $\co{2}$ in (\ref{eqn:embedding}).
On the other hand,
the embedding in (\ref{eqn:embedding}) appears 
in the construction of a unitary isomorphism
between the Bose-Fock space and the Fermi-Fock space
\cite{BFO01}.
From this, we wonder the existence of a behind mathematical structure 
among dynamical system, operator algebra
and quantum field theory.
\item
Apparently,
an embedding of $\coni$ into $\co{2}$ is not unique.
For example,
the following is one of other embeddings of $\coni$ into $\co{2}$:
\[s_{n}=t_{2}^{n-1}(t_{1}t_{2}t_{1}^{*}+t_{1}^{2}t_{2}^{*})\quad(n\geq 1).\]
Therefore the choice of the embedding (\ref{eqn:embedding}) 
is also an interesting accident.
\end{enumerate}
}
\end{rem}

In $\S$ \ref{section:second}, we prove Theorem \ref{Thm:maintwo}.
In $\S$ \ref{section:third}, we show examples of Theorem \ref{Thm:maintwo}.

%
%
\sftt{Proof of Theorem \ref{Thm:maintwo}}
\label{section:second}
We prove Theorem \ref{Thm:maintwo} in this section.
%
%
\begin{Thm}
\label{Thm:Sch}
(\cite{Schweiger}, $\S$ 19.2.1)
Let $(X,\mu,T)$ be a dynamical system and assume that
$\mu$ is $T$-invariant and
a measurable subset $A$ of $X$ satisfies (\ref{eqn:inverse}).
Define the new measure $\nu$ by
%
%
\begin{equation}
\label{eqn:nudash}
\nu(E)\equiv \mu(T^{-1}(E)\cap A)\quad(E\subset X).
\end{equation}
Then $\nu$ is invariant with respect to the jump transformation of $(T,A)$.
\end{Thm}
%
%
%
\begin{lem}
\label{lem:chain}
Let $(X,\mu,T)$ be a dynamical system,
and let $\nu_{1}$ and $\nu_{2}$ be other measures on $X$.
Assume that $\nu_{1}, \nu_{2}$ and $\mu\circ T$ are absolutely continuous 
with respect to $\mu$.
We write
$\psi_{i}\equiv d\nu_{i}/d\mu$ for $i=1,2$
and $\Phi_{T}\equiv d(\mu\circ T)/d\mu$.
Then the following are equivalent:
\begin{enumerate}
\item
$\nu_{1}=\nu_{2}\circ T$.
\item
$\psi_{1}=\Phi_{T}\cdot (\psi_{2}\circ T)$
almost everywhere in $X$.
\end{enumerate}
\end{lem}
%
%
\pr
By using the chain rule of \rnd s,
the statement holds.\qedh 

\noindent
From Lemma \ref{lem:chain}, the following holds.
%
%
\begin{cor}
\label{cor:function}
For $a,b\in {\bf R}$, $a<b$, let $[a,b]$ and $(a,b)$ denote the closed interval
and the open interval.
Let $T$ be a piecewise differentiable map from $(a,b)$ to $(a,b)$
with the differential $T^{'}$
and let $\nu_{1}$ and $\nu_{2}$ be measures on $[a,b]$.
Assume that $\nu_{1}, \nu_{2}$ and $\lambda\circ T$ are absolutely continuous 
with respect to the Lebesgue measure $\lambda$ on $[a,b]$.
Define
$\psi_{i}\equiv d\nu_{i}/d\lambda$ for $i=1,2$.
Then the following are equivalent:
\begin{enumerate}
\item
$\nu_{1}=\nu_{2}\circ T$.
\item
$\psi_{1}=|T^{'}|\cdot (\psi_{2}\circ T)$
almost everywhere in $[a,b]$
where $|T^{'}|(x)\equiv |T^{'}(x)|$ for $x\in (a,b)$.
\end{enumerate}
\end{cor}
\quad \\

\noindent
{\it Proof of Theorem \ref{Thm:maintwo}.}
(i)
Define the branching function system $g=\{g_{n}:n\in {\bf N}\}$ on $X$ by
%
%
\begin{equation}
\label{eqn:gmap}
g_{n}\equiv f_{2}^{n-1}\circ f_{1}\quad(n\geq 1).
\end{equation}
Since $f$ is a \bfs,
we can verify that $g$ is a branching function system on $X$.
From (\ref{eqn:embedding}) and (\ref{eqn:homomorphism}), 
(\ref{eqn:restriction}) holds.
If $G$ denotes the coding map of $g$,
then
%
%
\begin{equation}
\label{eqn:gmaptwo}
G(x)
=g_{n}^{-1}(x)
=(f_{1}^{-1}\circ f_{2}^{-n+1})(x)
\quad \mbox{ when }x\in g_{n}(X),\,n\in {\bf N}.
\end{equation}
Remark $g_{n}(X)=f_{2}^{n-1}(A_{1})$ for each $n\in {\bf N}$.
On the other hand,
when $x\in g_{n}(X)$, $e(x)=n-1$.
From this, $J_{\sigma,A_{1}}$ satisfies that
%
%
\begin{equation}
\label{eqn:tautwotwo}
J_{\sigma,A_{1}}(x)=\sigma^{n}(x)=(f_{1}^{-1}\circ f_{2}^{-n+1})(x).
\end{equation}
From (\ref{eqn:jump}),
$G=J_{\sigma,A_{1}}$ holds almost everywhere in $X$.
Hence the statement holds.

\noindent
(ii)
Define measures $\nu_{\infty}$ and $\nu_{2}$ of $X$ by
%
%
\begin{equation}
\label{eqn:nu}
\nu_{\infty}(E)\equiv 
\int_{E}\psi(x)\,d\mu(x),\quad
\nu_{2}(E)\equiv 
\int_{E}\phi(x)\,d\mu(x)
\end{equation}
for each measurable subset $E$ of $(X,\mu)$.
From (\ref{eqn:density}) and the definition of $\pi_{f}$,
%
%
\begin{equation}
\label{eqn:nutwo}
\psi(x)=\Bigl\{\{S(f_{1})^{*}\sqrt{\phi}\}(x)\Bigr\}^{2}
=\Phi_{f_{1}}(x)\cdot \phi(f_{1}(x)).
\end{equation}
From Lemma \ref{lem:chain},
%
%
\begin{equation}
\label{eqn:nuthree}
\nu_{\infty}(E)
=(\nu_{2}\circ f_{1})(E)
=\nu_{2}(\sigma^{-1}(E)\cap A_{1}).
\end{equation}
From Theorem \ref{Thm:Sch},
the statement holds.
\qedh

%
%
\sftt{Examples}
\label{section:third}
We show examples of Theorem \ref{Thm:maintwo} in this section.
For convenience,
we show  formulae of representations as follows.
Assume that $X=[a,b]$
and $T$ is a piecewise differentiable map on $X$
which is the coding map of a branching function system 
$\{f_{i}:i\in\Omega_{N}\}$ on $X$
for $2\leq N\leq \infty$.
Then (\ref{eqn:pfotwo}) is given as follows:
%
%
\begin{eqnarray}
\label{eqn:dif}
\{\pi_{f}(s_{i})\phi\}(x)=&\chi_{R_{i}}(x)\sqrt{|T^{'}(x)|}\,\phi(T(x)),\\
\nonumber
\\
\label{eqn:diftwo}
\{\pi_{f}(s_{i}^{*})\phi\}(x)=&\sqrt{|f^{'}_{i}(x)|}\,\phi(f_{i}(x))
\end{eqnarray}
for $i\in\Omega_{N}$.

%
%
\ssft{Tent map}
\label{subsection:thirdone}
As an introductory example of jump transformation,
we show a tent map.
Define the map $\Lambda$ by
$\Lambda(x)\equiv 1-|2x-1|$ on $x\in [0,1]$.
The map $\Lambda$ is called a {\it tent map} \cite{Devaney}.
Let $A_{1}\equiv [0,1/2]$, and
let $f_{2}\equiv (\Lambda|_{[0,1/2]})^{-1}$
and $f_{1}\equiv (\Lambda|_{[1/2,1]})^{-1}$.
Then
%
%
\begin{equation}
\label{eqn:tenttwo}
f_{1}(x)=1-\frac{1}{2}x,\quad
f_{2}(x)=\frac{1}{2}x.
\end{equation}
Define $g_{n}\equiv f_{2}^{n-1}\circ f_{1}$ for $n\geq 1$.
Then
%
%
\begin{equation}
\label{eqn:tentthree}
g_{n}(x)=\frac{2-x}{2^{n}}.
\end{equation}
The coding map $G$ of $\{g_{n}:n\in {\bf N}\}$ is given by
%
%
\begin{equation}
\label{eqn:tentfour}
G(x)=2-2^{n}x\quad\mbox{ when }1/2^{n}< x\leq 1/2^{n-1}.
\end{equation}
Hence the jump transformation is given as follows:
%
%
\begin{equation}
\label{eqn:tentfive}
J_{\Lambda,A_{1}}(x)=G(x)=2-2^{\lfloor \log_{2}x^{-1}\rfloor+1}\cdot x\quad(x\in (0,1]).
\end{equation}

\noindent
\def\pict{
\begin{minipage}{.45\linewidth}
\includegraphics[width=\linewidth]{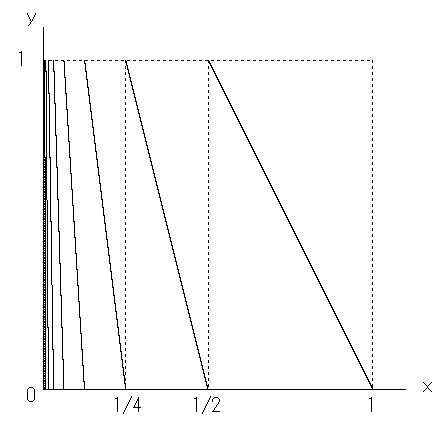}
\end{minipage}
}
\def\pictb{
\begin{minipage}{.45\linewidth}
\includegraphics[width=\linewidth]{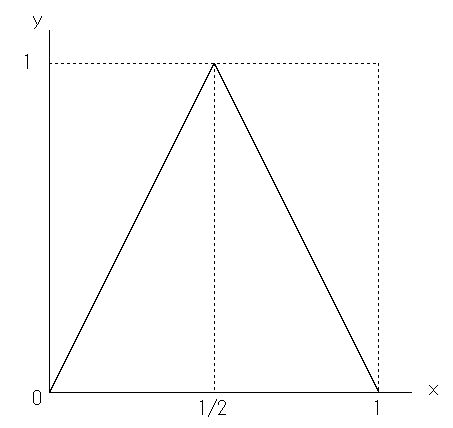}
\end{minipage}
}
\def\pictwo{
\put(0,350){\pict}
\put(200,600){$y=J_{\Lambda,A_{1}}(x)$}
}
\def\picthree{
\put(0,350){\pictb}
\put(200,600){$y=\Lambda(x)$}
}
\setlength{\unitlength}{.1mm}
\begin{picture}(1200,600)(0,50)
\put(600,0){\pictwo}
\put(0,0){\picthree}
\end{picture}

\noindent
Both of these  invariant probability measures are the Lebesgue measure.
For $f=\{f_{1},f_{2}\}$ in (\ref{eqn:tenttwo})
and $g=\{g_{n}:n\in {\bf N}\}$ in (\ref{eqn:tentthree}),
the representation $(L_{2}[0,1],\pi_{f})$ of $\co{2}$ 
and that $(L_{2}[0,1],\pi_{g})$ of $\coni$ are given as follows:
%
%
\begin{equation}
\label{eqn:tentsix}
\left\{
\begin{array}{rl}
\{\pi_{f}(t_{1})\phi\}(x)=&\chi_{[1/2,1]}(x)\sqrt{2}\phi(\Lambda(x)),\\
\\
\{\pi_{f}(t_{2})\phi\}(x)=&\chi_{[0,1/2]}(x)\sqrt{2}\phi(\Lambda(x)),\\
\end{array}
\right.
\end{equation}
%
%
\begin{equation}
\label{eqn:coni}
\{\pi_{g}(s_{n})\phi\}(x)=\chi_{[2^{-n},2^{-(n-1)}]}(x)2^{n/2}\phi(G(x))
\quad(n\geq 1).
\end{equation}
From (\ref{eqn:homomorphism}),
we can verify that
\[
\pi_{f}(s_{n})=
\pi_{f}(t_{2})^{n-1}\pi_{f}(t_{1})=
S(f_{2})^{n-1}S(f_{1})=
S(f_{2}^{n-1}\circ f_{1})=S(g_{n})=\pi_{g}(s_{n})
\]
for each $n$. Hence $\pi_{f}|_{\coni}=\pi_{g}$.

%
%
\ssft{Gauss map and Farey map}
\label{subsection:thirdtwo}
It is well known that the Gauss map $\tau_{1}$ in (\ref{eqn:gauss}) 
is the jump transformation of the {\it Farey map} $\sigma_{1}$ 
\cite{Thaler,KS02,Schweiger},  which is defined  on $[0,1]$ by
%
%
\begin{equation}
\label{eqn:tmaptwo}
\sigma_{1}(x)\equiv \frac{2}{1+|1-2x|}-1.
\end{equation}

\noindent
%
\def\pict{
\begin{minipage}{.4\linewidth}
\includegraphics[width=\linewidth]{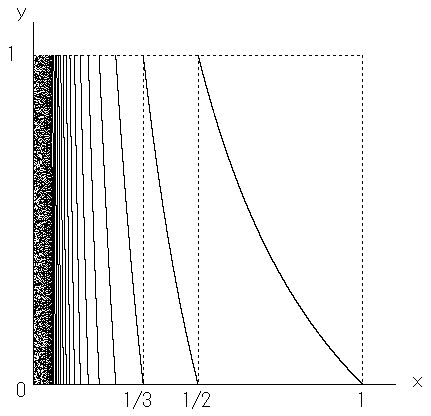}
\end{minipage}
}
\def\picttwo{
\begin{minipage}{.4\linewidth}
\includegraphics[width=\linewidth]{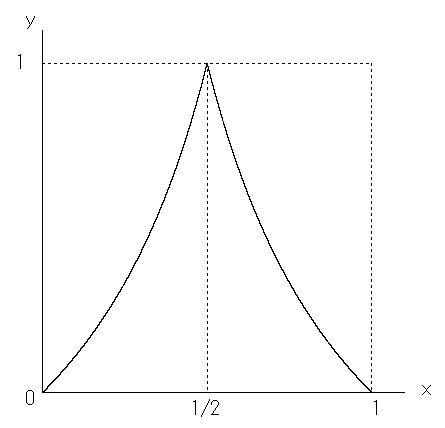}
\end{minipage}
}
\def\pictwo{
\put(50,350){\pict}
\put(-580,350){\picttwo}
\put(120,580){$y=\tau_{1}(x)=J_{\sigma_{1},A_{1}}(x)$}
\put(-420,580){$y=\sigma_{1}(x)$}
}
\setlength{\unitlength}{.1mm}
\begin{picture}(1200,600)(0,50)
\put(600,0){\pictwo}
\end{picture}

\noindent
where $A_{1}\equiv [0,1/2]$.

Define two measures $\gamma$ and $\theta$ on $[0,1]$ by
%
%
\begin{eqnarray}
\label{eqn:gamma}
d\gamma(x)\equiv &\disp{\frac{1}{\log 2}\frac{dx}{x+1},}\\
\nonumber
\\
\label{eqn:theta}
d\theta(x)\equiv &\disp{\frac{dx}{x}.}
\end{eqnarray}
The measure $\gamma$ is called {\it Gauss' measure} (\cite{IK}, p.16). 
Then $\gamma$ and $\theta$
are invariant  measures of $\tau_1$ and $\sigma_{1}$,
respectively  \cite{Thaler},
where $dx$ denotes the Lebesgue measure on $[0,1]$.
Remark that the former is finite but the later is not.

For $k\geq 1$,
define $g_{k}\equiv (\tau_{1})^{-1}$ on $\frac{1}{k+1}< x\leq \frac{1}{k}$,
and $f_{1}\equiv (\sigma_{1}|_{[1/2,1]})^{-1}$
and 
$f_{2}\equiv (\sigma_{1}|_{[0,1/2]})^{-1}$.
Then $\{g_{k}:k\in {\bf N}\}$ 
and $\{f_{1},f_{2}\}$ are
branching function systems
on $[0,1]$ with the coding map $\tau_{1}$ and $\sigma_{1}$,
respectively. 
Furthermore,
%
%
\begin{eqnarray}
\label{eqn:fk}
g_{k}(x)=\frac{1}{x+k}\quad(k\geq 1),\\
\nonumber
\\
\label{eqn:gtwo}
f_{1}(x)=\frac{1}{x+1},\quad
f_{2}(x)=\frac{x}{x+1}
\end{eqnarray}
for $x\in [0,1]$.
We see that $g_{k}=f_{2}^{k-1}\circ f_{1}$ for $k\geq 1$.
From these, we obtain representations $\pi_{\tau_{1}}$ 
and $\pi_{\sigma_{1}}$ of $\coni$ and $\co{2}$ on $L_{2}[0,1]$, respectively
which are written as follows:
%
%
\begin{equation}
\label{eqn:tauone}
\{\pi_{\tau_{1}}(s_{n})\phi\}(x)=
\frac{\chi_{[1/(n+1),1/n]}(x)}{x}\phi(\tau_{1}(x))\quad(n\in {\bf N}),
\end{equation}
%
%
\begin{equation}
\label{eqn:sigmaone}
\left\{
\begin{array}{rl}
\{\pi_{\sigma_{1}}(t_{1})\phi\}(x)=&
\disp{\frac{\chi_{[1/2,1]}(x)}{x}\phi(\sigma_{1}(x))},\\
\\
\{\pi_{\sigma_{1}}(t_{2})\phi\}(x)=&
\disp{\frac{\chi_{[0,1/2]}(x)}{1-x}\phi(\sigma_{1}(x))}
\end{array}
\right.
\end{equation}
for $x\in[0,1]$ and $\phi\in L_{2}[0,1]$.

%
%
\begin{rem}
\label{rem:gauss}
{\rm
Remark that 
the invariant measure $\theta$ of the Farey map in (\ref{eqn:theta}) is not finite.
Therefore Theorem \ref{Thm:maintwo}(ii) can not apply to this case.
However, if one forgets that the operator $\pi_{\sigma_{1}}(t_{1})^{*}$ 
is defined as an operator from $L_{2}[0,1]$ to $L_{2}[0,1]$,
then we see that Theorem \ref{Thm:maintwo}(ii) can  apply to also this case
as follows:
By definition,
%
%
\begin{equation}
\label{eqn:tstar}
\{\pi_{\sigma_{1}}(t_{1})^{*}\phi\}(x)=\frac{1}{x+1}\phi\left(\frac{1}{x+1}\right)
\quad(\phi\in L_{2}[0,1],\,x\in[0,1]).
\end{equation}
Although the density $\rho_{0}(x)\equiv \frac{1}{x}$  of $\theta$
does not belong to $L_{1}[0,1]$,
%
%
\begin{equation}
\label{eqn:tstartwo}
\Bigl\{\{\pi_{\sigma_{1}}(t_{1})^{*}\sqrt{\rho_{0}}\}(x)\Bigr\}^{2}=
\left\{
\frac{1}{x+1}
\sqrt{\frac{1}{\frac{1}{x+1}}}\right\}^{2}=\frac{1}{x+1}
\quad(x\in[0,1]).
\end{equation}
Therefore the
$\left\{\{\pi_{\sigma_{1}}(t_{1})^{*}\sqrt{\rho_{0}}\}(x)\right\}^{2}$
is the density of Gauss' measure $\gamma$ in (\ref{eqn:gamma}) 
up to scalar multiple.
Hence Theorem \ref{Thm:maintwo}(ii) holds for this case.
}
\end{rem}

We summarize relations between maps and Cuntz algebras
in this example as follows:\\

\begin{tabular}{l|l}
\hline
{\sf map} & {\sf Cuntz algebra}\\
\hline
Gauss map & representation of $\coni$\\
Farey map & representation of $\co{2}$\\
jump transformation & embedding of $\coni$ into $\co{2}$\\
\hline
\end{tabular}\\

%
%
\ssft{A generalization of continued fraction transformation}
\label{subsection:thirdthree}
In this subsection,
we show an example of jump transformation
associated with a generalization of continued fraction transformation.
Define the map $\tau_2$ on $[0,1]$ by
%
%
\begin{equation}
\label{eqn:ttwo}
\begin{array}{rl}
\tau_2(x)
=&\left\{
\begin{array}{ll}
\disp{\frac{1}{2^{k-1}x}-1}
\quad&(2^{-k}<x \leq 2^{-k+1},\,k\in {\bf N}),\\
\\
0\quad &(x=0).
\end{array}
\right.
\end{array}
\end{equation}
The map $\tau_{2}$ was introduced by Chan \cite{Chan}.
The invariant probability measure $\mu_{2}$ of 
$\tau_{2}$ is given as follows (\cite{Chan}, (2.20)):
%
%
\begin{equation}
\label{eqn:measureone}
d\mu_{2}(x)=\frac{1}{\log\frac{4}{3}}\frac{dx}{(x+1)(x+2)}\quad(x\in [0,1]).
\end{equation}
For $k\in {\bf N}$,
define $g_{k}\equiv (\tau_{2}|_{X_{k}})^{-1}$ on $X_{k}\equiv (2^{-k},2^{-k+1}]$.
Then
%
%
\begin{equation}
\label{eqn:inverseone}
g_{k}(x)=\frac{1}{2^{k-1}(x+1)}\quad(x\in [0,1)).
\end{equation}
Define the map $\sigma_{2}$ from $[0,1]$ to $[0,1]$ by
%
%
\begin{equation}
\label{eqn:sigmatwo}
\sigma_{2}(x)\equiv \left\{
\begin{array}{ll}
2x\quad &(0\leq x\leq \frac{1}{2}),\\
\\
\disp{\frac{1}{x}-1}\quad &(\frac{1}{2}\leq x\leq 1).
\end{array}
\right.
\end{equation}
Then $\tau_{2}$ is the jump transformation of $(\sigma_{2},A_{1})$ 
for $A_{1}\equiv [0,1/2]$, 
and Gauss' measure $\gamma$ in (\ref{eqn:gamma})
is $\sigma_{2}$-invariant.

\noindent
%
\def\pict{
\begin{minipage}{.4\linewidth}
\includegraphics[width=\linewidth]{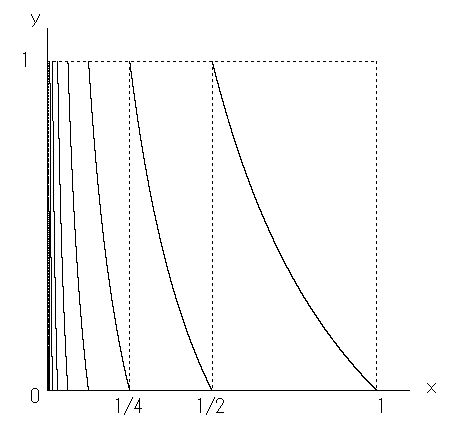}
\end{minipage}
}
\def\picttwo{
\begin{minipage}{.4\linewidth}
\includegraphics[width=\linewidth]{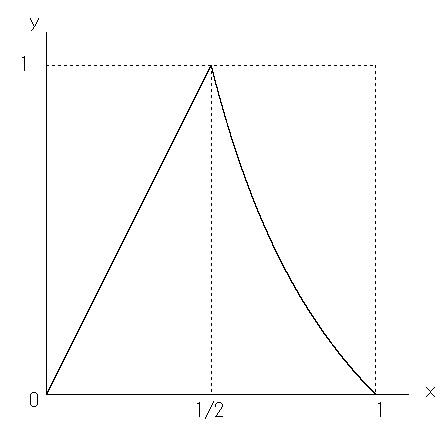}
\end{minipage}
}
\def\pictwo{
\put(50,350){\pict}
\put(-580,350){\picttwo}
\put(140,580){$y=\tau_{2}(x)=J_{\sigma_{2},A_{1}}(x)$}
\put(-410,580){$y=\sigma_{2}(x)$}
}
\setlength{\unitlength}{.1mm}
\begin{picture}(1200,600)(0,50)
\put(600,0){\pictwo}
\end{picture}

Define maps $f_{1}$ and $f_{2}$ on $[0,1]$ by 
%
%
\begin{equation}
\label{eqn:ftwo}
f_{1}(x)\equiv \frac{1}{x+1},\quad
f_{2}(x)\equiv \frac{1}{2}x.
\end{equation}
Then $\{f_{1},f_{2}\}$ is a branching function system
on $[0,1]$ with the coding map $\sigma_{2}$.
We see that
%
%
\begin{equation}
\label{eqn:mapf}
g_{k}(x)=f_{2}^{k-1}\circ f_{1}\quad(k\geq 1).
\end{equation}

For $f=\{f_{1},f_{2}\}$ in (\ref{eqn:ftwo}),
the representation $\pi_{f}$ of $\co{2}$
on $L_{2}[0,1]$ is given by
%
%
\begin{equation}
\label{eqn:rep}
\left\{
\begin{array}{rl}
\{\pi_{f}(t_{1})\phi\}(x)\equiv &
\chi_{[1/2,1]}(x)(1+x)\phi\left(\frac{1}{x}-1\right),
\\
\\
\{\pi_{f}(t_{2})\phi\}(x)\equiv &\chi_{[0,1/2]}(x)\sqrt{2}\phi(2x).
\end{array}
\right.
\end{equation}
for $\phi \in L_{2}[0,1]$ and $x\in [0,1]$.
From this,
%
%
\begin{equation}
\label{eqn:reptwo}
\left\{
\begin{array}{rl}
\{\pi_{f}(t_{1}^{*})\phi\}(x)=& \frac{1}{1+x}\phi\left(\frac{1}{x+1}\right),\\
\\
\{\pi_{f}(t_{2}^{*})\phi\}(x)= &\frac{1}{\sqrt{2}}\phi\left(\frac{1}{2}x\right).
\end{array}
\right.
\end{equation}

Define $\phi_{0}(x)\equiv \frac{1}{\log 2}\frac{1}{x+1}$.
From Theorem \ref{Thm:maintwo}(ii),
\[\psi(x)=
\Bigl\{\{\pi(t_{1}^{*})\sqrt{\phi_{0}}\}(x)\Bigr\}^{2}
=\frac{1}{(1+x)^{2}}\frac{1}{1+\frac{1}{1+x}}
=\frac{1}{(1+x)(2+x)}.\]
In this way,
the the density $\psi$ of $\tau_{2}$-invariant measure $\mu_{2}$ in (\ref{eqn:measureone})
is given up to scalar multiple.

For further generalizations of $\tau_{1}$ and $\tau_{2}$,
see \cite{CL}.
About operator theory associated with continued fractions,
see \cite{IK,Mayer,Schweiger}.


%


%
\end{document}